# BOUNDED SOLUTIONS TO BACKWARD SDE'S WITH JUMPS FOR UTILITY OPTIMIZATION AND INDIFFERENCE HEDGING

By Dirk Becherer

*Imperial College, London*

We prove results on bounded solutions to backward stochastic equations driven by random measures. Those bounded BSDE solutions are then applied to solve different stochastic optimization problems with exponential utility in models where the underlying filtration is noncontinuous. This includes results on portfolio optimization under an additional liability and on dynamic utility indifference valuation and partial hedging in incomplete financial markets which are exposed to risk from unpredictable events. In particular, we characterize the limiting behavior of the utility indifference hedging strategy and of the indifference value process for vanishing risk aversion.

**1. Introduction.** A prominent stochastic control problem in stochastic finance is the utility maximization problem, where the objective is to maximize by optimal investment the expected utility from future wealth. Another problem is the valuation and hedging of contingent claims in incomplete markets. Here the task is to determine jointly a suitable notion of dynamic valuation and an optimal partial hedging strategy such that both are consistent with no-arbitrage theory. The utility indifference approach combines the two aforementioned problems and has recently attracted a lot of interest. Because of many recent publications on this approach, we refrain from giving yet another survey on the topic and the literature, but refer for introductions with more comprehensive references to the articles by Fittelli [11], Delbaen et al. [8], Becherer [2] and Mania and Schweizer [21]. In this article, we obtain solutions for the exponential utility maximization problem under an additional liability, both on the primal and on the dual level, and for the utility indifference valuation and hedging problem in a financial market model that permits for nontradable risk from unpredictable events. This event risk can involve both the nonpredictable time of an event and a possibly nonpredictable event size. An example is provided by the jump times and sizes of a marked point process. In such a model, we obtain limiting results for vanishing risk aversion for both the dynamic utility indifference value process and for the corresponding indifference hedging strategy. Typical areas of application may be models from the areas of credit risk or from insurance, where doubly marked point processes are commonly used to model the occurrence and size of default or insurance losses.

Our main mathematical means to solve these optimization problems are results on backward stochastic equations (BSDEs) with jumps, which are derived in the first part of the paper. BSDEs are generally known to be useful for studying problems in mathematical finance (see [9]), but have been mainly used in continuous settings thus far. Broadly, BSDEs can describe optimality equations from dynamical programming; this is similar to the familiar Hamilton–Jacobi–Bellman equations, but more general in that BSDEs can also cover non-Markovian situations. The second part of the paper shows how the solutions to our two stochastic optimization problems with exponential utility can be described explicitly in terms of our BSDE solutions. Our incomplete market framework does not necessitate the underlying filtration to be continuous, but allows for noncontinuous martingales and nonpredictable stopping times; this motivates the first contribution of the paper on BSDEs that are driven by a Brownian motion and a random measure. Building on results of Tang and Li [24] and Barles, Buckdahn and Pardoux [1], we derive existence, uniqueness and continuity results for bounded solutions to such BSDEs when the generator possesses a certain monotonicity. This extends previous results on square integrable solutions to solutions with more integrability and with a possibly nonhomogeneous random measure, making them amendable to the subsequent applications.

Our second contribution is the solution, by an application of our BSDE results, of the two aforementioned exponential utility optimization problems in a model with nonpredictable jump risk. The articles by Rouge and El Karoui [22], Hu, Imkeller and Müller [14] and Mania and Schweizer [21] are closely related to this part of of paper. While their framework basically covers continuous filtrations, the prime example being the Brownian filtration, the present paper works in a setting with random measures which allows the modeling of risk from nonpredicable events. We moreover characterize the limit of the indifference hedging strategy for vanishing risk aversion; this complements a recent result by Mania and Schweizer [21] for continuous filtrations, and our analysis shows how a random measure component in the BSDE provides the natural means to solve similar problems in the presence of unpredictable jump risk. Another contribution is the method of proof which derives the solution to our optimization problems directly from our existence and uniqueness results for bounded BSDEs with random measures.

Having proven existence and uniqueness of a suitable BSDE solution, we can employ the martingale optimality principle to solve the (primal) optimal stochastic control problems directly. This is more in the spirit of the work by Rouge and El Karoui [22] and Hu, Imkeller and Müller [14], whereas Mania and Schweizer [21] derive existence and uniqueness for their specific BSDEs from the given existence of optimal solutions, ensured by duality results of Kabanov and Stricker [16]. By *continuity results* for BSDEs with jumps, we obtain the asymptotic behavior of the solution to the indifference valuation and hedging problem for vanishing risk aversion. The limit corresponds to risk minimization under the minimal entropy martingale measure.

The paper is organized as follows. Section 2 establishes the general framework, assuming the existence of a stochastic basis carrying a Brownian motion and a compensated integer-valued random measure that possess a weak predictable representation property. Section 3 derives existence, uniqueness and continuity results for bounded solutions of BSDEs with jumps whose generator may not satisfy the usual global Lipschitz condition. Section 4 applies these results to study the exponential utility maximization problem with an additional liability and the utility indifference valuation and hedging problem.

**2. Framework and preliminaries.** This section sets out the notation and the assumptions that are supposed to hold in the sequel.

We start with a stochastic basis $(\Omega, \mathcal{F}, \mathbb{F}, P)$ with a finite time horizon $T < \infty$ and a filtration $\mathbb{F} = (\mathcal{F}_t)_{t \in [0,T]}$ satisfying the usual conditions of right continuity and completeness, such that we can and do take all semimartingales to have right continuous paths with left limits. For simplicity, we assume that $\mathcal{F}_0$ is trivial and $\mathcal{F} = \mathcal{F}_T$. Conditional expectations with respect to $\mathcal{F}_t$ (and $P$) are denoted by $E_t[\cdot] = E_t^P[\cdot]$. On this stochastic basis, let $W = (W_t)$ be a $d$-dimensional standard Brownian motion and let $\mu$ denote an integer-valued random measure

$$\mu(dt, de) = (\mu(\omega, dt, de) | \omega \in \Omega)$$

on $([0,T] \times E, \mathcal{B}([0,T]) \otimes \mathcal{E})$ with compensator

$$\nu := \nu^P(dt, de)$$

under $P$, where $E := \mathbb{R}^\ell \setminus \{0\}$ is equipped with its Borel $\sigma$-field $\mathcal{E} := \mathcal{B}(E)$. Define the measure $P \otimes \nu$ on $(\widetilde{\Omega}, \widetilde{\mathcal{F}}) := (\Omega \times [0,T] \times E, \mathcal{F} \otimes \mathcal{B}([0,T]) \otimes \mathcal{E})$ by

$$(2.1) \qquad P \otimes \nu(\widetilde{B}) = E\left[\int_{[0,T] \times E} I_{\widetilde{B}}(\omega, t, e) \nu(\omega, dt, de)\right], \qquad \widetilde{B} \in \widetilde{\mathcal{F}};$$

this is called the *measure generated by* $\nu$. Let $\mathcal{P}$ denote the predicable $\sigma$-field on $\Omega \times [0,T]$ and define

$$\widetilde{\mathcal{P}} := \mathcal{P} \otimes \mathcal{E}.$$

A function on $\widetilde{\Omega}$ that is $\widetilde{\mathcal{P}}$-measurable is called *predictable*. We suppose that $\nu$ is equivalent to a product measure $\lambda \otimes dt$ with a density $\zeta$ such that

$$(2.2) \qquad \nu(\omega, dt, de) = \zeta(\omega, t, e)\lambda(de)\, dt,$$

where $\lambda$ is a $\sigma$-finite measure on $(E, \mathcal{E})$ satisfying $\int_E 1 \wedge |e|^2 \lambda(de) < \infty$, and where the density $\zeta$ is a $\widetilde{\mathcal{P}}$-measurable, bounded, nonnegative function such that for some, constant $c_\nu$,

$$(2.3) \qquad 0 \leq \zeta(\omega, t, e) \leq c_\nu < \infty, \qquad P \otimes \lambda \otimes dt\text{-a.e.}$$

Clearly, (2.2) implies that $\nu(\{t\} \times E) = 0$ for all $t$, and $\nu([0,T] \times E) \leq c_\nu T \lambda(E)$.

For a predicable function $U$ on $\widetilde{\Omega}$, the integral process with respect to $\mu$ (analogously for $\nu$) is defined as

$$U * \mu_t(\omega) = \begin{cases} \int_{[0,t] \times E} U(\omega, s, e)\mu(\omega, ds, de), & \text{if finitely defined}, \\ +\infty & \text{otherwise}. \end{cases}$$

We recall that for any predictable function $U$, the process $U * \nu$ is a predictable process, while the process $U * \mu$ is an optional process, and that $E[|U| * \mu_T] = E[|U| * \nu_T]$. If $(|U|^2 * \mu)^{1/2}$ is locally integrable, then $U$ is integrable with respect to $\widetilde{\mu} = \mu - \nu$, and $U \stackrel{(P)}{*} \widetilde{\mu}$ is defined as the purely discontinuous local martingale (under $P$) with jump process $(\int_E U \mu(\{t\}, de))_t$. If, moreover, the process $|U|^2 * \nu$ is integrable, then $U$ is integrable with respect to $\widetilde{\mu}$, and $U * \widetilde{\mu} = U * (\mu - \nu)$ is a square integrable, purely discontinuous martingale with predictable quadratic variation $\langle U * (\mu - \nu) \rangle = |U|^2 * \nu$. If the increasing process $|U| * \mu$ (or, equivalently $|U| * \nu$) is locally integrable, then $U$ is $\widetilde{\mu}$-integrable and $U * \widetilde{\mu} = U * \mu - U * \nu$. We refer to [15] for details on (integer-valued) random measures and stochastic integrals and note that our assumptions on $\mu$ and $\nu$ imply that $\widehat{W} = 0$ in Section II.1.d of [15].

We assume that, with respect to $\mathbb{F}$ and $P$,

(2.4)  $W$ and $\widetilde{\mu}$ have the weak property of predictable representation.

This means that every square integrable martingale $M$ has a representation

$$(2.5) \qquad M = M_0 + Z \cdot W + U * \widetilde{\mu} := M_0 + \int Z\, dW + U * \widetilde{\mu},$$

where $Z$ and $U : \widetilde{\Omega} \to \mathbb{R}$ are predictable processes such that $E[\int_0^T |Z|^2\, dt] < \infty$ and $E[|U|^2 * \nu_T] < \infty$, that is, both stochastic integrals are in the space $\mathcal{H}^2$ of square integrable martingales. We next provide several cases of interest where (2.4) holds:

EXAMPLE 2.1. (1) Let $W$ be a Brownian motion and let $N$ be an independent Poisson point process. Then $W$ and the compensated measure $\widetilde{\mu}^N$ of the jump measure $\mu^N$ of $N$, that is $\mu^N(dt, de) := \sum_{s \in (0,T]} \delta_{(s, \Delta N_s)}(dt, de) I_{\{\Delta N_s \neq 0\}}$, have the representation property (2.4) with respect to the usual filtration $\mathbb{F}^{(W,N)}$ generated by them.

(2) More generally, let $(X_t)_{t \in [0,T]}$ be a marked point process, that is, a process whose paths are RCLL step functions with only a finite number of jumps (i.e., $X$ can be represented as $X = x_0 + \sum_i \xi_i \mathbb{1}_{[\![T_i, T]\!]}$ with $x_0 \in \mathbb{R}^\ell$ and with random times $T_i \in (0, \infty]$ such that $T_i \uparrow \infty$ and $0 < T_i < T_{i+1}$ on $\{T_i < \infty\}$ for all $i$, where all marks $\xi_i$ are $\mathbb{R}^\ell$-valued random variables with $\{\xi_i = 0\} = \{T_i = \infty\}$). Let $W$ be a Brownian motion independent of $X$ and let $\mathbb{F} := \mathbb{F}^{(W,X)}$ denote the usual filtration generated by $X$ and $W$. Let $\mu := \mu^X := \sum_i \delta_{(T_i, \xi_i)}(dt, de) I_{\{T_i \leq T\}}$ and let $\nu$ denote the compensator, being the same under $\mathbb{F}^X$ and $\mathbb{F}$. Then $W$ and $\widetilde{\mu}$ have property (2.4) with respect to $\mathbb{F}$. To see this, note that both $W$ and $\widetilde{\mu}$ have the representation property with respect to their own filtrations (see [13], Theorems 13.19, 5.52). By strong orthogonality, each martingale $M$ with $M_T = I_A I_B$ for $A \in \mathcal{F}_T^W$ and $B \in \mathcal{F}_T^X$ can be represented as in (2.5). This implies (2.4), since the linear span of random variables like $\mathbb{1}_A \mathbb{1}_B$ is dense in $L^2(\mathcal{F}_T)$.

(3) Let $X$ be a (time-homogenous) Lévy process with $X_0 = 0$ and predictable characteristics $(\alpha, \beta, \nu)$. Then the continuous martingale part $X^c$ and the compensated jump measure $\widetilde{\mu}^X = \mu^X - \nu$ of $X$ have the representation property (2.5) with respect to the filtration $\mathbb{F}^X$; see [13], Theorems 13.44 and 13.49. If $X^c$ does not vanish (i.e, $\beta \neq 0$), then there is a constant $C \in (0, \infty)$ such that $W := X^c / C$ is a Brownian motion.

(4) Suppose $W$ and $\widetilde{\mu}$ have the representation property (2.4) under $P$. Let $P'$ be a probability measure absolutely continuous (or equivalent) to $P$ with density process $(Z_t)_{t \in [0,T]}$. Then the $P'$-Brownian motion $W' := W - \int (Z_-)^{-1} d\langle Z, W \rangle$ and the $P'$-compensated jump measure $\widetilde{\mu}' = \mu^X - \nu^{P'}$ have the representation property (2.4) with respect to $(\Omega, \mathbb{F}, P')$; see Theorem 13.22 in [13]. This offers plenty of scope to build models where $W$ and $\nu$ are not independent from the previous examples.

**3. Backward stochastic differential equations with jumps.** For ease of exposition, all results in this section are formulated with $P$ representing some generic probability measure.

REMARK 3.1. The results of the present section will be used in the sequel with different equivalent measures taking the role of $P$. This causes no problems when these measures, the corresponding Brownian motions and the compensators for $\mu$ satisfy the same assumptions as imposed on $P$, $W$ and $\nu$ (cf. Example 2.1 4). It will be made clear on those occasions with respect to which measure the results and notation are to be used.

Let us fix some notation:

- $\mathcal{S}^{p,k}$ with $1 \leq p \leq \infty$ denotes the space of $\mathbb{R}^k$-valued semimartingales $(Y_t)_{t \in [0,T]}$ with $\|Y\|_{\mathcal{S}^p} := \|\sup_{t \in [0,T]} |Y_t|\|_{L^p} < \infty$.
- $\mathcal{L}_T^{2,k}$ ($\mathcal{L}_T^{2,k \times d}$) denotes the space of $\mathcal{P}$-predictable processes $Z$ taking values in $\mathbb{R}^k$ ($\mathbb{R}^{k \times d}$) with $\|Z\|_{\mathcal{L}_T^2} := (E[\int_0^T |Z_t|^2 \, dt])^{1/2} < \infty$. This norm is equivalent to the norm $\|Z\|_\beta := \|(e^{\beta t} Z_t)_{t \in [0,T]}\|_{\mathcal{L}_T^2}$ for $\beta \in \mathbb{R}$.
- $\mathcal{L}_\nu^{2,k}$ denotes the space of $\widetilde{\mathcal{P}}$-predictable functions $U : \widetilde{\Omega} \to \mathbb{R}^k$ with

$$\|U\|_{\mathcal{L}_\nu^2} := \left( E\left[ \int_0^T \int_E |U_t(e)|^2 \nu(dt, de) \right] \right)^{1/2} < \infty.$$

Assumptions (2.2) and (2.3) imply that the space $\mathcal{L}_\nu^{2,k} = L^2(\widetilde{\mathcal{P}}, P \otimes \nu; \mathbb{R}^k)$ includes the space $\mathcal{L}_{\lambda \times dt}^{2,k} = L^2(\widetilde{\mathcal{P}}, P \otimes \lambda \otimes dt; \mathbb{R}^k)$.

- $L^0(\mathcal{E}, \lambda; \mathbb{R}^k)$ denotes the space of measurable functions with the topology of convergence in measure. It will be convenient to define for $u, u' \in L^0(\mathcal{E}, \lambda : \mathbb{R}^k)$,

$$(3.1) \qquad \|u - u'\|_t := \left( \int_E |u - u'|^2 \zeta(t, e) \lambda(de) \right)^{1/2}.$$

For $U \in \mathcal{L}_\nu^{2,k}$, $\|U_t\|_t < \infty$ holds $P \otimes dt$-a.e. as $E[\int \|U_t\|_t^2 \, dt] = \|U\|_{\mathcal{L}_\nu^2}^2$.

To simplify notation we will omit dimension indices like $k$ when they are clear from the context. But in later sections we shall refer in our notation to the underlying probability measure when it is different from $P$.

For a given data tuple $(B, f)$, which consists of a random variable $B$ and a suitable generator function $f_t(y, z, u) = f(\omega, t, y, z, u)$, we are interested in finding a triple $(Y, Z, U)$ of processes in a suitable space such that

$$(3.2) \quad Y_T = B \quad \text{and} \quad dY_t = -f_t(Y_{t-}, Z_t, U_t) \, dt + Z_t \, dW_t + \int_E U_t(e) \widetilde{\mu}(dt, de)$$

for $t \in [0, T]$. Equation (3.2) can be written in integrated form as

$$(3.3) \quad Y_t = B + \int_t^T f_s(Y_{s-}, Z_s, U_s) \, ds - \int_t^T Z_s \, dW_s - \int_t^T \int_E U_s(e) \widetilde{\mu}(ds, de),$$

with $t \in [0, T]$. Such a triple $(Y, Z, U)$ is called a *solution* to the backward stochastic differential equation (3.2) or (3.3).

Proposition 3.2 ensures existence and uniqueness of the BSDE solution in an $L^2$-sense for our setting with a nonhomogeneous compensator $\nu$. Admitting $\nu$ to be nonhomogeneous allows more interesting mutual dependencies to occur between the tradable and nontradable risk factors in our later applications. For the homogeneous case with $\zeta \equiv 1$ in (2.2), the result was shown in Lemma 2.4 of [24]; see also Theorem 2.1 in [1]. It is straightforward to generalize the established fixed point method of proof to the present setting, hence we leave details to the reader.

PROPOSITION 3.2. *Let $B \in L^2(\mathcal{F}_T, P; \mathbb{R}^k)$ and suppose that the function*
$$f : \Omega \times [0,T] \times \mathbb{R}^k \times \mathbb{R}^{k \times d} \times L^0(\mathcal{E}, \lambda; \mathbb{R}^k) \to \bar{\mathbb{R}}^k$$
*is $\mathcal{P} \otimes \mathcal{B}(\mathbb{R}^k) \otimes \mathcal{B}(\mathbb{R}^{k \times d}) \otimes \mathcal{B}(L^0(\mathcal{E}, \lambda; \mathbb{R}^k))$-measurable, satisfies $f_t(0,0,0) \in \mathcal{L}_T^{2,k}$ and that there exists a constant $K_f \in [0, \infty)$ such that*

(3.4) $\quad |f_t(y,z,u) - f_t(y',z',u')| \leq K_f(|y-y'| + |z-z'| + \|u-u'\|_t)$

*holds $P \otimes dt$-a.e. for all $y, y' \in \mathbb{R}^k$, $z, z' \in \mathbb{R}^{k \times d}$, and $u, u' \in L^0(\mathcal{E}, \lambda; \mathbb{R}^k)$ [in particular, the left-hand side of (3.4) is supposed to be finite when the right-hand side is finite]. Then there exists a unique $(Y, Z, U)$ in $\mathcal{S}^{2,k} \times \mathcal{L}_T^{2,k \times d} \times \mathcal{L}_\nu^{2,k}$ which solves the BSDE (3.2).*

The next continuity result generalizes Proposition 2.2 from [1] to a conditional estimate, which reduces to the unconditional estimate from [1] for $\tau = 0$. Let us remark that $f$ could take nonfinite values for some $u \in L^0$, but (3.4) implies that for all $(Y, Z, U)$ in $\mathcal{S}^{2,k} \times \mathcal{L}_T^{2,k \times d} \times \mathcal{L}_\nu^{2,k}$, the process $(f(Y_-, Z, U))$ is finite $P \otimes dt$-a.e. Hence, the $\delta f$ term in (3.5) is well defined $P \otimes dt$-a.e.

PROPOSITION 3.3. *Let $(B, f)$ and $(B', f')$ be data satisfying the assumptions of Proposition 3.2, with solutions $(Y, Z, U)$ and $(Y', Z', U')$ in $\mathcal{S}^{2,k} \times \mathcal{L}_T^{2,k \times d} \times \mathcal{L}_\nu^{2,k}$, respectively. Denote $(\delta B, \delta f) = (B - B', f - f')$ and let $(\delta Y, \delta Z, \delta U) = (Y - Y', Z - Z', U - U')$. Then there exists a constant $c = c(T, K_{f'}) < \infty$ depending on $T$ and $K_{f'}$ such that*

$$(3.5) \quad \begin{aligned} E_\tau &\left[ \sup_{u \in [\![\tau,T]\!]} |\delta Y_u|^2 + \int_\tau^T |\delta Z_s|^2 \, ds + \int_{]\!]\tau,T]\!] \times E} |\delta U_s|^2 \nu(ds, de) \right] \\ &\leq c E_\tau \left[ |\delta B|^2 + \int_\tau^T |\delta f_s(Y_{s-}, Z_s, U_s)|^2 \, ds \right] < \infty \end{aligned}$$

*holds for all stopping times $\tau \leq T$. If, moreover, the random variable $\delta B$ and the process $(\delta f(Y_-, Z, U))$ are bounded, then $\int \delta Z \, dW$ and $\delta U * \widetilde{\mu}$ are $\mathrm{BMO}(P)$-martingales.*

Noting that the BSDE has a trivial solution for vanishing data, one sees that estimate (3.5) with $\tau = t$ implies the useful a priori estimate

$$(3.6) \quad \begin{aligned} E_t &\left[ \sup_{t \leq u \leq T} |Y_u|^2 + \int_t^T |Z_s|^2 \, ds + \int_{(t,T] \times E} |U_s|^2 \nu(ds, de) \right] \\ &\leq c E_t \left[ |B|^2 + \int_t^T |f_s(0,0,0)|^2 \, ds \right] < \infty. \end{aligned}$$

PROOF OF PROPOSITION 3.3. We extend the argument from [1]. Applying Itô's formula to $|\delta Y|^2$ yields that, for $t \in [0, T]$,

$$|\delta Y_t|^2 - |\delta B|^2 + \int_t^T |\delta Z_s|^2 + \|\delta U_s\|_s^2 \, ds$$

$$(3.7) \quad - 2 \int_t^T \delta Y_{s-}(f_s(Y_{s-}, Z_s, U_s) - f'(Y'_{s-}, Z'_s, U'_s)) \, ds$$

$$= -2 \int_t^T \delta Y_{s-} \delta Z_s \, dW_s - \int_t^T \int_E 2\delta Y_{s-} \delta U_s + (\delta U_s)^2 \widetilde{\mu}(ds, de).$$

Using the integrability of $(\delta Y, \delta Z, \delta U) \in \mathcal{S}^{2,k} \times \mathcal{L}_T^{2,k \times d} \times \mathcal{L}_\nu^{2,k}$ and the assumptions imposed on $f$ and $f'$, the left-hand side of equation (3.7) can be dominated in absolute value by an integrable random variable, uniformly in $t$. Hence, the stochastic integrals on the right-hand side are in the martingale space $\mathcal{H}^1$ and their increments vanish in conditional expectation. Using this fact, the inequality $2Kab \leq 2K^2a^2 + b^2/2$ ($K, a, b \in \mathbb{R}$) and the Lipschitz condition on $f'$ yields, for any stopping time $\tau \leq T$ and $u \in [0, T]$, that

$$E_\tau \left[ |\delta Y_{\tau \vee u}|^2 + \int_{\tau \vee u}^T |\delta Z_s|^2 + \|\delta U_s\|_s^2 \, ds \right]$$

$$= E_\tau \left[ |\delta B|^2 + \int_{\tau \vee u}^T |\delta f_s(Y_{s-}, Z_s, U_s)|^2 \, ds. \right.$$

$$\left. + c \int_{\tau \vee u}^T |\delta Y_s|^2 \, ds + \tfrac{1}{2} \int_{\tau \vee u}^T |\delta Z_s|^2 + \|\delta U_s\|_s^2 \, ds \right]$$

with $\tau \vee u := \max(\tau, u)$ and $c$ denoting a constant depending on $T$ and $K_{f'}$. By the conditional Fubini result of Lemma A.1 (with $\sigma = \tau$) and Gronwall's lemma, at this point we nearly obtain the desired inequality (3.5), but with a $\sup_{u \in [0,T]}$ outside the conditional expectation. From the BSDE that $\delta Y$ satisfies, we have

$$|\delta Y_{\tau \vee u}| \leq E_{\tau \vee u} \left[ |\delta B| + \int_\tau^T |f_s(Y_{s-}, Z_s, U_s) - f'_s(Y'_{s-}, Z'_s, U'_s)| \, ds \right],$$

for $u \in [0, T]$. Taking the supremum over $u \in [0, T]$ on both sides and applying Doob's inequality to the supremum of the $(\mathcal{F}_{\tau \vee u})_{u \in [0, T]}$-martingale on the right-hand side yields that $E_\tau[\sup_{u \in [\![\tau, T]\!]} |\delta Y_u|^2]$ is dominated by the term $cE_\tau[|\delta B|^2 + \int_\tau^T |\delta f_s(Y_{s-}, Z_s, U_s)|^2 + |f'_s(Y_{s-}, Z_s, U_s) - f'_s(Y'_{s-}, Z'_s, U'_s)|^2 \, ds]$ with some constant $c$ depending on $T$. Estimate (3.5) now follows from the Lipschitz property of $f'$ and the previous estimate.

Finally, if $\delta B$ and $(\delta f_t(Y_-, Z, U))_t$ are bounded, then it follows from its BSDE that $\delta Y$ is bounded, and that the left-hand side of (3.5) is bounded by a constant, uniformly in $\tau$. Hence, the integral processes $\int \delta Z \, dW$ and $\delta U * \widetilde{\mu}$ are *BMO*-martingales; see Theorem 10.9.4 in [13]. □

For future reference, we state the following simple but useful result:

LEMMA 3.4. *Let $Y \in \mathcal{S}^\infty$ and $B \in L^\infty$. Suppose that $Z$ is a $\mathcal{P}$-predictable process and $U$ is a $\widetilde{\mathcal{P}}$-predictable function which are integrable in the sense of local martingales with respect to $W$ and $\widetilde{\mu}$, respectively, and $f(\omega, t, y, z, u)$ is a product-measurable function such that $f(Y_-, Z, U)$ is in $L^\infty(P \otimes dt)$ and that $(Y, Z, U)$ solves the BSDE (3.2) with data $(B, f)$. Then the stochastic integrals $\int Z\, dW$ and $U * \widetilde{\mu}$ are both BMO(P)-martingales. In particular, $Z \in \mathcal{L}_T^2$ and $U \in \mathcal{L}_\nu^2$.*

PROOF. Being a bounded martingale, the process $\int Z\, dW + U * \widetilde{\mu}$ is in BMO($P$). By the characterization of BMO-martingales (see [13], Theorem 10.9), the claim then follows from (i) the observation that the quadratic covariation of $\int Z\, dW$ and $U * \widetilde{\mu}$ vanishes, so that the quadratic variation of each addend is dominated by the quadratic variation of the sum, and (ii) from the fact that the jumps of $U * \widetilde{\mu}$ are the jumps of $Y$, hence bounded. □

In our later applications of BSDEs, the Lipschitz condition (3.4) on the generator $f$ will not be satisfied, and more than square integrability of $Y$ will be needed. To this end, we show that, basically, a monotonicity property (3.11) of the generator with respect to the jumps, together with bounded terminal data, ensures existence of a bounded solution to the BSDE (3.2). Let us emphasize that for BSDEs with jumps, a comparison result does not, in general, hold under the assumptions of Proposition 3.2; see the example in [1]. Therefore, we can not infer the existence of a bounded solution for bounded terminal data by such means. For BSDEs without jumps that are driven solely by Brownian motions and which have quadratic generators that may not be globally Lipschitz, existence results were obtained in [18].

For the remainder, we consider the one-dimensional case with $k = 1$ and generator functions $f : \Omega \times [0, T] \times \mathbb{R}^1 \times \mathbb{R}^{1 \times d} \times L^0(\mathcal{E}, \lambda; \mathbb{R}) \to \bar{\mathbb{R}}$ of the form

$$f_t(y, z, u) = \begin{cases} \widehat{f}_t(y, z, u) + \int_E g_t(u(e))\zeta(t, e)\lambda(de), & \text{if finitely defined,} \\ +\infty, & \text{otherwise,} \end{cases}$$

(3.8)

where $\widehat{f}$ and $g$ satisfy conditions which, although $f$ does not, in general, satisfy the assumptions of Proposition 3.2, still ensure existence and uniqueness of a BSDE solution whose components $Y$ and $U$ are furthermore bounded.

THEOREM 3.5. *Let $k = 1$. Assume that $B \in L^\infty$ is bounded, $\lambda(E) < \infty$ is finite, and that $f$ has the form (3.8) where $\widehat{f}$ satisfies all assumptions from*

*Proposition* 3.2 *(for f). Assume additionally that there exist* $K_1, K_2 \in [0, \infty)$ *such that*

(3.9) $\quad |\widehat{f}_t(y, z, u)| \leq K_1 + K_2 |y| \qquad$ *holds* $P \otimes dt$-*a.e. for all* $y, z, u$,

*and that* $g : \Omega \times [0, T] \times \mathbb{R} \to \mathbb{R}$ *in (3.8) is a* $\mathcal{P} \otimes \mathcal{B}(\mathbb{R})$-*measurable function such that*

(3.10) $\quad u \mapsto g_t(u) \quad$ *is locally Lipschitz in* $u$, *uniformly in* $(\omega, t)$, *and*

(3.11) $\begin{aligned} g_t(u) &\leq -u = +|u| \text{ for } u \leq 0, \\ g_t(u) &\geq -u = -|u| \text{ for } u \geq 0, \end{aligned} \qquad P \otimes dt$-*a.e.*

*Then there exists a unique solution* $(Y, Z, U)$ *in* $\mathcal{S}^\infty \times \mathcal{L}^2_T \times \mathcal{L}^2_\nu$ *to the BSDE (3.2), where* $Y \in \mathcal{S}^\infty$ *is bounded and* $U$ *is bounded* $P \otimes \nu$-*a.e.*

*Moreover, if* $(B', f')$ *is another tuple of data satisfying the assumption of this theorem with solution* $(Y', Z', U') \in \mathcal{S}^\infty \times \mathcal{L}^2_T \times \mathcal{L}^2_\nu$, *then estimate (3.5) of Proposition* 3.3 *still holds. In particular, estimate (3.6) still holds.*

Under the assumptions of this theorem, one can thus choose representatives for $Y$ and $U$ which are bounded on $\Omega \times [0, T]$ and $\widetilde{\Omega}$, respectively. Before proving the theorem, let us give an example of a generator function which will reappear in the later applications to mathematical finance.

EXAMPLE 3.6. For $\alpha \in (0, \infty)$, let $g(\omega, t, u) = g_t(u) := \frac{1}{\alpha} e^{\alpha u} - u - \frac{1}{\alpha}$ and $\widehat{f} := 0$. Then $f$ from (3.8) satisfies the assumptions of Theorem 3.5.

PROOF OF THEOREM 3.5. By hypothesis, there exists $K_3 \in [0, \infty)$ such that $|B| \leq K_3$. Define a truncation-boundary function $b : [0, T] \to \mathbb{R}^+$ by

(3.12) $\quad b(t) := \begin{cases} K_3 + K_1(T - t), & \text{when } K_2 = 0, \\ K_3 e^{K_2(T-t)} + \dfrac{K_1}{K_2}(e^{K_2(T-t)} - 1), & \text{when } K_2 > 0, \end{cases}$

and a truncation function

(3.13) $\qquad \kappa(t, y) := \min(\max(y, -b(t)), +b(t)).$

Then define $\widetilde{f} : \Omega \times [0, T] \times \mathbb{R} \times \mathbb{R}^d \times L^0(\mathcal{E}, \lambda; \mathbb{R}^1) \to \bar{\mathbb{R}}$ by

$$\widetilde{f}(\omega, t, y, z, u) = \widehat{f}(\omega, t, \kappa(t, y), z, \kappa(t, y + u) - \kappa(t, y))$$
$$+ \int_E g_t(\kappa(t, y + u(e)) - \kappa(t, y)) \zeta(t, e) \lambda(de)$$

when the right-hand side is finite, and by $+\infty$ elsewhere. Using (3.13), the local Lipschitz property of $g$ and the fact that $\lambda(E) < \infty$, one can verify by the Schwartz inequality that $\widetilde{f}$ satisfies the assumptions of Proposition 3.2.

We note that one could find examples of $\lambda$, if $\lambda(E)$ were not finite, where the Lipschitz condition (3.4) would not be met by $\widetilde{f}$ for $f$ from Example 3.6. Let $(Y, Z, U) \in \mathcal{S}^2 \times \mathcal{L}_T^2 \times \mathcal{L}_\nu^2$ denote the unique solution for the BSDE with generator $\widetilde{f}$. For $\widetilde{Y}_t := \kappa(t, Y_t)$, let

$$\widetilde{U}_t(e) := \kappa(t, Y_{t-} + U_t(e)) - \kappa(t, Y_{t-}),$$

which represents the jump size of $\widetilde{Y}$. We will show below that

(3.14)
$$\text{the processes } Y \text{ and } \widetilde{Y} \text{ are indistinguishable and}$$
$$U = \widetilde{U} \text{ holds } P \otimes \nu\text{-a.e.}$$

By definition, $|\widetilde{Y}_t| \leq b(t) \leq b(0)$ and $|\widetilde{U}_t(e)| \leq 2b(t) \leq 2b(0)$ are bounded uniformly in $t$. Clearly, (3.14) implies that $Y = \widetilde{Y}$ in $\mathcal{S}^2 \subset \mathcal{L}_T^2$ and $U = \widetilde{U}$ in $\mathcal{L}_\nu^2$, hence $\int \widetilde{f}_t(Y_{t-}, Z_t, U_t) \, dt = \int \widetilde{f}_t(\widetilde{Y}_{t-}, Z_t, \widetilde{U}_t) \, dt$ and $U * \widetilde{\mu} = \widetilde{U} * \widetilde{\mu}$.

Admitting result (3.14) for a moment, it follows that the solution $(Y, Z, U)$ of the BSDE with generator $\widetilde{f}$ also solves the BSDE with $f$, and that $Y$ is in $\mathcal{S}^\infty$ and $U$ has a bounded representative $\widetilde{U}$ in $\mathcal{L}_\nu^2$. To show uniqueness, let $(Y', Z', U')$ be a another solution to the BSDE with $f$, with $Y'$ bounded. Similarly as with $\widetilde{U}$ for $U$, one can find a bounded representative for $U'$ in $\mathcal{L}_\nu^2$, for example, bounded by $2\|Y'\|_{\mathcal{S}^\infty}$, using the fact that $\int_E |U'_t(e)| \mu(\{t\}, de) = |\Delta Y'_t| \leq 2\|Y'\|_{\mathcal{S}^\infty}$. By taking $K_3$ larger when necessary, say $K_3 \geq 2\|Y'\|_{\mathcal{S}^\infty}$, one can assume that $|Y'|, |U'| \leq K_3$. But then both $(Y, Z, U)$ and $(Y', Z', U')$ solve the BSDE also with generator $\widetilde{f}$, and by the uniqueness result from Proposition 3.2 applied to $\widetilde{f}$, the two solutions must coincide. Finally, the validity of estimate (3.6) and the validity of the estimate from Proposition 3.3 follow from the observation that the BSDE solutions to $(f, B)$ and $(f', B')$ also solve the BSDEs with the corresponding truncated generators $\widetilde{f}$ and $\widetilde{f}'$.

To complete the proof, it remains to show (3.14). To prove $|Y_t| \leq b(t)$ for all $t$, we first consider the upper bound. Fix $t \in [0, T]$ and let

$$\tau := \inf\{s \in [t, T] | Y_s \leq b(s)\}.$$

Then $Y_s \geq b(s)$ for $(\omega, s) \in [\![t, \tau [\![$, and one has $Y_\tau \leq b(\tau)$ with $\tau \leq T$, since $Y_T \leq K_3 = b(T)$. Since $(Y, Z, U)$ solves the BSDE with $\widetilde{f}$, it follows that

$$Y_t = E_t \left[ Y_\tau + \int_t^\tau \widetilde{f}(s, Y_{s-}, Z_s, U_s) \, ds \right]$$
$$\leq E_t \left[ Y_\tau + \int_t^\tau K_1 + K_2 b(s) \, ds \right.$$
$$\left. + \int_t^\tau \int_E g_s(\kappa(s, Y_{s-} + U(e)) - \kappa(s, Y_{s-})) \mu(ds, de) \right]$$

$$\leq E_t\left[Y_\tau + \int_t^\tau K_1 + K_2 b(s)\,ds + g_\tau(\kappa(\tau,Y_\tau) - \kappa(\tau,Y_{\tau-}))\right]$$

$$\leq E_t\left[Y_\tau + \int_t^\tau (K_1 + K_2 b(s))\,ds + b(\tau) - Y_\tau\right]$$

$$\leq E_t\left[\int_t^\tau (K_1 + K_2 b(s))\,ds + b(\tau)\right] = b(t);$$

this uses the fact that $G * \widetilde{\mu} = G * \mu - G * \nu$ is a martingale for

$$G_s := g_s(\kappa(s, Y_{s-} + U(e)) - \kappa(s, Y_{s-})),$$

since $|G_s| * \nu \leq c(|U_s| \wedge 2b(s)) * \nu$ is integrable (see [15], II.1.28–30), where $c < \infty$ denotes a local Lipschitz constant for $g$. We have further used the fact that

$$\int_{]\!]t,\tau]\!]\times E} G_s \mu(ds, de) = \sum_{s \in ]\!]t,\tau]\!]} g_s(\kappa(s, Y_s) - \kappa(s, Y_{s-}))\mathbb{1}_{(\Delta Y_s \neq 0)},$$

with all summands vanishing on $s \in ]\!]t,\tau[\![$ except possibly at $s = \tau$, where

$$g_\tau(\kappa(\tau, Y_\tau) - \kappa(\tau, Y_{\tau-})) \leq \kappa(\tau, Y_{\tau-}) - \kappa(\tau, Y_\tau) \leq b(\tau) - Y_\tau.$$

This shows the upper bound $Y_t \leq b(t)$ for any $t$. The lower bound $Y_t \geq -b(t)$ is proved similarly by using $\tau := \inf\{s \in [t, T] | Y_s \geq -b(s)\}$. Since $Y$ has RCLL-paths, it follows that $Y$ and $\widetilde{Y}$ are indistinguishable processes. Hence,

$$0 = \sum_{t \in (0,T]} (\Delta(Y_t - \widetilde{Y}_t))^2 = \int_0^T \int_E (U(e) - \widetilde{U}(e))^2 \mu(dt, de) = (U - \widetilde{U})^2 * \mu_T,$$

implying that $E[(U - \widetilde{U})^2 * \mu_T] = E[(U - \widetilde{U})^2 * \nu_T]$ vanishes. This establishes the second part of (3.14).  □

**4. Applications in exponential utility optimization.** This section applies the previous results to solve two prominent optimization problems with exponential utility. We start with the expected utility maximization problem with an additional liability, and proceed afterwards to the utility indifference valuation and hedging problem.

4.1. *The financial market framework.* Within the general framework of Section 2, we now introduce a financial market model. In this section, the measure $P$ represents the objective "true world" probability measure. All assumptions and notation of the present subsection remain valid for the remainder of Section 4, in addition to those of the general framework.

The market contains a riskless numeraire asset as well as $d$ risky assets, whose discounted price processes $S = (S_t^i)_{t \in [0,T]}$, $i = 1, \ldots, d$, evolve in $(0, \infty)^d$ according to the stochastic differential equation

$$
\begin{aligned}
(4.1) \quad & dS_t = \operatorname{diag}(S_t^i)_{i=1,\ldots,d} \sigma_t (\varphi_t \, dt + dW_t) =: \Sigma_t \, d\widehat{W}_t, \qquad t \in [0, T], \\
& S_0 \in (0, \infty)^d,
\end{aligned}
$$

where $\varphi$ is a predictable from $\mathbb{R}^d$-valued $dt$-integrable process, and $\sigma$ is an $\mathbb{R}^{d \times d}$-valued predictable process such that $\sigma_t$ is invertible $P \otimes dt$-a.e. and integrable with respect to

$$
(4.2) \qquad \widehat{W} := W + \int \varphi_t \, dt,
$$

and we define $\Sigma_t := (\operatorname{diag}(S_t^i)_{i=1,\ldots,d}) \sigma_t$. We suppose that

(4.3)     the market price of risk process $\varphi$ is bounded $P \otimes dt$-a.e.

The solution $S$ to (4.1) is given by the stochastic exponential

$$
S_t^i = S_0^i \mathcal{E}\left(\int \sum_{j=1}^d \sigma_t^{ij} (\varphi^j \, dt + dW_t^j)\right)_t = S_0^i \mathcal{E}\left(\left(\int \sigma_t \, d\widehat{W}_t\right)^i\right)_t, \qquad t \in [0, T].
$$

The notation $\widehat{W}$ is reminiscent of the fact that $\widehat{W}$ is a Brownian motion with respect to the so-called minimal martingale measure

$$
(4.4) \qquad d\widehat{P} := \mathcal{E}\left(-\int \varphi_t \, dW\right)_T dP.
$$

By (4.4), it follows that

(4.5)   the compensator of the random measure $\mu$ under $\widehat{P}$ equals $\nu$,

that is, it equals the compensator under $P$, up to indistinguishability. We assume

$$
(4.6) \qquad\qquad\qquad \lambda(E) < \infty,
$$

such that $\nu([0, T] \times E)$ is bounded by (2.2) and (2.3).

For any $S$-integrable $\mathbb{R}^d$-valued process $\vartheta$, the gains process from trading according to a strategy to hold $\vartheta$ shares of risky assets $S$ is given by $\int \vartheta \, dS = \int \theta \, d\widehat{W}$, where

$$
(4.7) \qquad\qquad \theta(\vartheta) := \Sigma^{\mathrm{tr}} \vartheta \quad \text{and} \quad \vartheta(\theta) = (\Sigma^{\mathrm{tr}})^{-1} \theta
$$

provides a bijection between $\vartheta$ and $\theta$. Because the parameterization with respect to $\theta$ eases notation and facilitates later formulas, we will use relation (4.7) to parameterize our strategies in terms of $\theta$ in the sequel.

In Section 4.2 below, we are going to consider an investor who wants to maximize the exponential utility with risk aversion $\alpha > 0$ from his terminal wealth. For such an investor, we define the set of available trading strategies

$$\Theta := \Theta(P, \alpha) \tag{4.8}$$

as follows. Let $\Theta$ consist of all $\mathbb{R}^d$-valued, predictable, $S$-integrable processes $\theta$ which meet the following integrability requirements under $P$:

$$E^P\left[\int_0^T |\theta_t|^2 dt\right] < \infty \tag{4.9}$$

and

$$(4.10) \quad \left\{\exp\left(-\alpha \int_0^\tau \theta_t d\widehat{W}_t\right)\Big|\tau \in \mathfrak{T}\right\} \text{ is a uniformly } P\text{-integrable family}$$

of random variables, with $\mathfrak{T}$ denoting the set of all stopping times $\tau \leq T$. Conditions (4.9) and (4.10) correspond to those in [14]; see their Definition 1 and subsequent remarks. Condition (4.9) excludes arbitrage possibilities like doubling strategies from $\Theta$, as explained below, and the exponential condition (4.10) fits rather naturally with the exponential preferences of our investor. We will also show that it transforms in a "good" way under a change to the minimal entropy martingale measure; see part (3) of Remark 4.5.

Let $\mathbb{P}_e := \{Q \sim P | S \text{ is a local } Q\text{-martingale}\}$ denote the set of all equivalent local martingale measures. The market is free of arbitrage in the sense that there exists at least one measure in $\mathbb{P}_e$, namely $\widehat{P}$, that has finite relative entropy with respect to $P$, thus

$$\mathbb{P}_f := \{Q \in \mathbb{P}_e | H(Q|P) < \infty\} \neq \varnothing. \tag{4.11}$$

Also, the set $\Theta$, whose strategies can give rise to wealth processes unbounded from below, does not contain arbitrage strategies, since $\int_0^T \theta\, d\widehat{W} \geq 0$ with $\theta \in \Theta$ implies that $\int \theta\, d\widehat{W} = 0$. In fact, $d\widehat{P}/dP$ and $(\int_0^T |\theta|^2 dt)^{1/2}$ are in $L^2(P)$ by (4.3), (4.4) and (4.9). So, Hölder's inequality yields that $(\int_0^T |\theta|^2 dt)^{1/2}$ is in $L^1(\widehat{P})$, hence $\int \theta\, d\widehat{W}$ is in the martingale space $\mathcal{H}^1(\widehat{P})$. The claim thus follows.

The market model is in general incomplete because we only assume the weak representation property (2.4) and do not assume a continuous or Brownian filtration. Despite the invertibility of the volatility matrix $\sigma$, there will in general be several martingale measures having different compensators for the random measure $\mu$, which represents some nontradable risk factors related to nonpredictable events such as the jump times and sizes of a marked (doubly) Poisson process. If $\mu$ is nontrivial, one can also see directly that there exist purely discontinuous martingales under any $Q$ in $\mathbb{P}_e$ which cannot be represented as stochastic integrals with respect to the continuous process $S$.

As contingent claims, we consider European claims payable at time $T$ whose payoff is described by a bounded random variable $B \in L^\infty := L^\infty(P)$. This integrability assumption is the same as in the articles [22], [14] and [21] which are closely related to the subsequent analysis, and it fits comfortably into the duality setting from [8].

4.2. *Exponential utility maximization.* We are going to consider the problem of maximizing the expected utility from terminal wealth at time $T$ for the exponential utility function $x \mapsto -\exp(-\alpha x)$ with risk aversion parameter $\alpha \in (0, \infty)$. We will later compare results for different levels of risk aversion and, to this end, we emphasize that most quantities in the sequel depend on $\alpha$. But notational references to $\alpha$ are omitted where $\alpha$ is clear from the context, as in most proofs.

Recalling our parameterization (4.7) of strategies, the solution of the problem of maximizing the expected exponential utility from terminal wealth

$$V_t^{B,\alpha}(x) := V_t^B(x)$$

$$(4.12) \qquad := \operatorname*{ess\,sup}_{\theta \in \Theta(P,\alpha)} E_t^P\left[-\exp\left(-\alpha\left(x + \int_t^T \vartheta(\theta)\,dS - B\right)\right)\right]$$

$$= \exp(-\alpha x) \operatorname*{ess\,sup}_{\theta \in \Theta(P,\alpha)} E_t^P\left[-\exp\left(-\alpha\left(\int_t^T \theta\,d\widehat{W} - B\right)\right)\right]$$

from ($\mathcal{F}_t$-measurable) capital $x \in \mathbb{R}$ at time $t \in [0,T]$ by optimal future investments under an additional liability $B$ will be described by the following BSDE under $P$:

$$Y_t = B + \int_t^T -Z_s\varphi_s - \frac{|\varphi|^2}{2\alpha}\,ds$$

$$(4.13) \qquad + \int_t^T \int_E \left(\frac{\exp(\alpha U_s(e))-1}{\alpha} - U_s(e)\right)\zeta(s,e)\lambda(de)\,ds$$

$$- \int_t^T Z_s\,dW_s - \int_t^T \int_E U_s(e)\widetilde{\mu}(ds,de), \qquad t \in [0,T].$$

The linear dependence of the generator in (4.13) on $Z$ can be removed by a change of measure from $P$ to $\widehat{P}$ from (4.4). By Theorem 12.29 of [13], the integral with respect to $\widetilde{\mu}$ remains unaltered by this change of measure because of the unchanged compensator (4.5). Hence the BSDE (4.13) transforms to the following BSDE under the measure $\widehat{P}$:

$$Y_t = B + \int_t^T -\frac{|\varphi|^2}{2\alpha}\,ds$$

$$(4.14) \qquad + \int_t^T \int_E \left(\frac{\exp(\alpha U_s(e))-1}{\alpha} - U_s(e)\right)\zeta(s,e)\lambda(de)\,ds$$

$$-\int_t^T Z_s \, d\widehat{W}_s - \int_t^T \!\!\int_E U_s(e)\widetilde{\mu}(ds, de), \qquad t \in [0, T].$$

THEOREM 4.1. *The solution to the utility maximization problem (4.12) with risk aversion $\alpha > 0$ and liability $B \in L^\infty$ is described by the unique solution*

$$(Y^B, Z^B, U^B) := (Y^{B,\alpha}, Z^{B,\alpha}, U^{B,\alpha}) \in \mathcal{S}^\infty(\widehat{P}) \times \mathcal{L}^2_T(\widehat{P}) \times \mathcal{L}^2_\nu(\widehat{P})$$

*to the BSDE (4.14) under $\widehat{P}$ [solving (4.13) under $P$, with $(Y^B, Z^B, U^B)$ being in $\mathcal{S}^\infty(P) \times \mathcal{L}^2_T(P) \times \mathcal{L}^2_\nu(P)$]. The optimal value function $V^{B,\alpha}_t(x)$ and the optimal strategy are given by*

(4.15) $\qquad V^{B,\alpha}_t(x) = -e^{-\alpha x} \exp(\alpha Y^B_t) = -e^{-\alpha x}\exp(-\alpha(-Y^B_t)),$

(4.16) $\qquad \theta^B := \theta^{B,\alpha} := Z^B + \dfrac{1}{\alpha}\varphi \in \Theta(P, \alpha).$

*Moreover, $\int \theta^B_t \, dW_t$ is in $BMO(P)$.*

In the sense of (4.15), $-Y^B$ can be considered as the exponential time-$t$-certain equivalent wealth, and $Y^B$ as the time-$t$-certain liability, which are equivalent to the gains (and losses) arising from both the future optimal investment and from the terminal liability.

PROOF OF THEOREM 4.1. Let us note that $(\widehat{W}, \nu)$ under $\widehat{P}$ fits into the setting of Section 3. The proof then proceeds in several steps.

First, Theorem 3.5 ensures the existence of the unique (bounded) solution $(Y, Z, U)$ in $\mathcal{S}^\infty(\widehat{P}) \times \mathcal{L}^2_T(\widehat{P}) \times \mathcal{L}^2_\nu(\widehat{P})$ to the BSDE (4.14) under $\widehat{P}$, with $U$ being bounded $P \otimes \nu$-a.e. and $\int Z \, d\widehat{W}$ being in $BMO(\widehat{P})$, by Lemma 3.4. By Theorem 3.6 in [17] and (4.4) it follows that $\int Z \, dW \in BMO(P) \subset \mathcal{H}^2(P)$ and $Z \in \mathcal{L}^2_T(P)$. Hence, $(Y, Z, U)$ also solves BSDE (4.13) under $P$, with the integrability claims for $Y$ and $U$ under $P$ following from their a.s. boundedness.

By the Doléans–Dade formula, or direct computation, one obtains

$$\mathcal{E}\left(\int\!\!\int_E \exp(\alpha U_s(e)) - 1 \widetilde{\mu}(ds, de)\right)_t$$
$$= \exp\left(\int_0^t \!\!\int_E \alpha U_s(e)\widetilde{\mu}(ds, de) \right.$$
$$\left. - \int_0^t \!\!\int_E \exp(\alpha U_s(e)) - 1 - \alpha U_s(e)\nu(ds, de)\right).$$

By the BSDE (4.14), equality (4.2) and by Itô's formula it follows that

$$
\begin{aligned}
(4.17) \quad &-\exp\left(-\alpha\left(Y_0 + \int_0^t \theta\, d\widehat{W} - Y_t\right)\right) \\
&= -\exp\left(-\alpha \int_0^t \theta - Z\, d\widehat{W} + \tfrac{1}{2}\int_0^t |\varphi|^2\, dt\right) \\
&\quad \times \mathcal{E}\left(\int\!\!\int \exp(\alpha U_s(e)) - 1\, \widetilde{\mu}(ds, de)\right)_t \\
&= -e^{(\alpha^2/2)\int_0^t |\theta - Z - \varphi/\alpha|^2\, dt}\, \mathcal{E}\left(-\alpha \int \theta - Z\, dW\right)_t \\
&\quad \times \mathcal{E}\left(\int\!\!\int \exp(\alpha U_s(e)) - 1\, \widetilde{\mu}(ds, de)\right)_t \\
&= -e^{(\alpha^2/2)\int_0^t |\theta - Z - \varphi/\alpha|^2\, dt} \\
(4.18) \quad &\quad \times \mathcal{E}\left(-\alpha \int \theta - Z\, dW + \int\!\!\int \exp(\alpha U_s(e)) - 1\, \widetilde{\mu}(ds, de)\right)_t,
\end{aligned}
$$

for any $\theta \in \Theta$, $t \in [0, T]$, where (4.18) uses Yor's formula. For any $\theta \in \Theta$,

(4.19) the stochastic exponential in (4.18) ($t \in [0, T]$) is a P-martingale.

To see the latter, note that condition (4.10) in the definition of $\Theta$ and the boundedness of $Y$ imply that the term in (4.17) is a uniformly integrable family for $t \in [0, T]$. Since the ordinary exponential factor in (4.18) is monotone and clearly bounded away from zero, the stochastic exponential in (4.18) is a uniformly integrable local martingale, hence (4.19).

By (4.18) and (4.19), the process from (4.18) is a supermartingale for all $\theta \in \Theta$ and a martingale for $\theta^B = Z + \varphi/\alpha$. This implies optimality of $\theta^B$, provided that we can show $\theta^B$ is in $\Theta$. To this end, observe that $\int \theta^B\, d\widehat{W}$ is in $BMO(\widehat{P})$, since $\int Z\, d\widehat{W}$ and $\int \varphi\, d\widehat{W}$ are by Lemma 3.4 and (4.3). Since $dP/d\widehat{P} \in L^p(\widehat{P})$ for any $1 \leq p < \infty$, $\int_0^T |\theta^B|^2\, dt$ is in $L^1(P)$ by Hölder's inequality. This yields (4.9), while (4.10) follows from (4.19) and the fact that $Y \in \mathcal{S}^\infty$. Moreover, $\int \theta^B\, dW$ is in $BMO(P)$ by Theorem 3.6 in [17].

Taking the conditional expectation of the utility of the optimal terminal wealth and using the martingale property of the process (4.18) for $\theta = \theta^B$ finally yields (4.15). $\square$

We have solved the primal utility maximization problem directly by using the classical (super-)martingale verification argument to show optimality of the candidate solution. To link our results to the martingale duality results, let us recall from [8] (cf. [2] or [21]) that the primal exponential utility maximization problem with liability $B$ is related to the dual problem of finding, for a given $\alpha > 0$,

$$(4.20) \qquad Q^{E,B} = \underset{Q \in \mathbb{P}_f}{\arg\max} \{\alpha E^Q[B] - H(Q|P)\}.$$

For $B = 0$, this means finding the minimal entropy martingale measure $Q^E := Q^{E,0}$, while $Q^{E,B}$ can be shown to minimize the relative entropy $H(Q|P_B)$ with respect to $dP_B := \mathrm{const}\, \exp(\alpha B)\, dP$ over the set $\mathbb{P}_f$.

The next theorem describes the density process for the solution $Q^{E,B}$ to the dual problem as an ordinary and also as a stochastic exponential, explicit in terms of the ingredients of the related BSDE. Furthermore, it describes the compensator of $\nu$ under $Q^{E,B}$. The density process of $Q^{E,B}$ turns out to be

$$
\begin{aligned}
Z_t^{E,B} &:= \exp\left(-\alpha\left(Y_0^B + \int_0^t \theta^B\, d\widehat{W} - Y_t^B\right)\right) \\
(4.21) \qquad &= \mathcal{E}\left(-\int \varphi\, dW + \iint_E \exp(\alpha U_s^B(e)) - 1\, \widetilde{\mu}(ds, de)\right)_t \\
&= \mathcal{E}\left(-\int \varphi\, dW\right)_t \mathcal{E}\left(\iint_E \exp(\alpha U_s^B(e)) - 1\, \widetilde{\mu}(ds, de)\right)_t,
\end{aligned}
$$

with $Y^B$, $U^B$ and $\theta^B$ from Theorem 4.1; the equalities hold by (4.17) and (4.18) for $\theta = \theta^B$. We note that even the ordinary exponential form does not follow immediately from existing duality results in [8] and elsewhere, since our definition of the strategy set $\Theta$ involves only integrability assumptions under the objective measure $P$; this differs from the variants of $\Theta$ studied in [8].

THEOREM 4.2. *Suppose the assumptions of Theorem* 4.1 *hold. Then the density process with respect to $P$ of the solution $Q^{E,B} \in \mathbb{P}_f$ to the dual problem (4.20) is given by $Z^{E,B}$ from (4.21). Furthermore, $\widehat{W}$ is a $Q^{E,B}$-Brownian motion, and the compensator of $\mu$ under $Q^{E,B}$ is given by*

$$(4.22) \qquad \nu^{Q^{E,B}}(dt, de) = \exp(\alpha U^B(e)) \nu(dt, de).$$

PROOF. We apply the verification theorem from [12] to identify the solution to the dual problem. To this end, we must validate that $Z_T^{E,B}$ satisfies three conditions. First, it is clear from Theorem 4.1 and (4.19) that $Z^{E,B}$ is a strictly positive density process, and the BSDE for $Y^B$ implies that

$$(4.23) \qquad Z_T^{E,B} = \exp\left(-\alpha\left(Y_0^B + \int_0^T \theta^B\, d\widehat{W} - B\right)\right).$$

Define

$$d\bar{Q} := \bar{Z}_T\, dP \qquad \text{for } \bar{Z} := Z^{E,B}.$$

By the stochastic exponential form of the density process from (4.21) and Girsanov's theorem, it follows that $\widehat{W}$ is a $\bar{Q}$-Brownian motion and that

$$\bar{\nu}(dt, de) := \exp(\alpha U^B(e))\nu(dt, de) \tag{4.24}$$

is the $\bar{Q}$-compensator of $\mu$. The first claim is standard. For the second claim, let $w$ denote an $\widetilde{\mathcal{P}}$-predictable function on $\widetilde{\Omega}$ such that $|w| * \mu$ is locally $\bar{Q}$-integrable. The latter is equivalent to $|w| * \mu$ being locally $P$-integrable, since $\bar{Z}$ and $1/\bar{Z}$ are both locally bounded due to the boundedness of $U$. By the form (4.21) of $\bar{Z}$ and the BSDE for $Y^B$, it follows that $\Delta\bar{Z}_t$ equals

$$\int_E \bar{Z}_{t-}(\exp(\alpha U^B(e)_t) - 1)\mu(\{t\}, de) = \Delta((\bar{Z}_-(\exp(\alpha U^B) - 1)) * \mu)_t.$$

By integration by parts and Propositions II.1.28 and 30 in [15], it follows that

$$d(\bar{Z}(w * \mu))_t = d(\bar{Z}_- w) * \mu_t + (w * \mu_{t-}) d\bar{Z}_t + \Delta\bar{Z}_t \Delta(w * \mu_t)$$
$$= d(\bar{Z}_- \exp(\alpha U^B)w) * \mu_t + (w * \mu_{t-}) d\bar{Z}_t,$$
$$d(\bar{Z}(w * \bar{\nu}))_t = d(\bar{Z}((w\exp(\alpha U^B)) * \nu))_t$$
$$= d(\bar{Z}_- \exp(\alpha U^B)w) * \nu_t + ((w\exp(\alpha U^B)) * \nu_{t-}) d\bar{Z}_t + 0.$$

By subtracting the two processes and using the fact that $\nu$ is the $P$-compensator of $\mu$, one obtains that $w * \mu - w * \bar{\nu}$ is a local $\bar{Q}$-martingale. Hence, $\bar{\nu}$ is the compensator of $\mu$ under $\bar{Q}$ by Theorem II.1.8 of [15].

Denoting $\bar{\mu} = \mu - \bar{\nu}$ and recalling the compensator relation (4.24), we have (see Theorem 12.28 from [13]) that under a change of measure,

$$U \overset{\widehat{P}}{*} \widetilde{\mu} = U \overset{\bar{Q}}{*} \bar{\mu} + (U(\exp(\alpha U) - 1)) * \nu, \tag{4.25}$$

with $U$ being integrable with respect to $\bar{\mu}$ such that $U * \bar{\mu}$ is a local martingale under $\bar{Q}$. Hence, it follows from the BSDE (4.14) under $\widehat{P}$ that $(Y^B, Z^B, U^B)$ is also a bounded solution to the following BSDE under $\bar{Q}$:

$$\begin{aligned}Y_t = B - \int_t^T Z_s \, d\widehat{W}_s - \int_t^T \int_E U_s(e)\bar{\mu}(ds, de) + \int_t^T \frac{-|\varphi|^2}{2\alpha} \, ds \\ + \int_t^T \int_E \left(\frac{\exp(\alpha U_s(e)) - 1}{\alpha} - U_s(e)\exp(\alpha U_s(e))\right)\zeta(s, e)\lambda(de) \, ds.\end{aligned} \tag{4.26}$$

Lemma 3.4 ensures integrability of $(Y^B, Z^B, U^B)$ under the change of measure. By (4.16), it follows that $\int \theta^B \, d\widehat{W}$ is a $\bar{Q}$-BMO-martingale, since

$$\int Z^B \, d\widehat{W} \quad \text{is a } \bar{Q}\text{-BMO-martingale} \tag{4.27}$$

and $\varphi$ is bounded. This is the second condition needed. For the third, note that Lemma 3.4 and the BSDE (4.14) under $\widehat{P}$ imply that $\int Z^B \, d\widehat{W}$, and hence $\int \theta^B \, d\widehat{W}$, is in $\mathrm{BMO}(\widehat{P})$. By the John–Nirenberg inequality, we conclude that there is some $\varepsilon > 0$ such that $\exp(\varepsilon \int_0^T \theta^B \, d\widehat{W})$ is in $L^p(\widehat{P})$ for some $p > 1$. Since $\phi$ is bounded, $dP/d\widehat{P} = \mathcal{E}(\int \phi \, d\widehat{W})_T$ is in $L^q(\widehat{P})$ for any $q \in [1, \infty)$, so

$$(4.28) \qquad \exp\left(\varepsilon \int_0^T \theta^B \, d\widehat{W}\right) \quad \text{is in } L^1(P)$$

by Hölder's inequality. By (4.23), (4.27) and (4.28), all three conditions for the verification result from [12] (cf. Proposition 3.5 in [2]) are satisfied, implying that $\bar{Q}$ is the optimal measure $Q^{E,B}$ and $\bar{Z} = Z^{E,B}$ is its density process. $\square$

REMARK 4.3. Let us point out some connections to the pioneering work in [6, 7]. In Sections 1.3 and 12, the authors study a model with jump risk that can be accommodated in our general framework [see Example 2.1(2) and (4)] it involves a single risky asset price following a geometric Brownian motion and one additional unpredictable (default) event with an intensity. They showed how the solution to the exponential utility problem can be derived from a certain BSDE under special, partially restrictive, assumptions, and noted in particular that a rigorous general existence result for a (sufficiently nice) solution to the key BSDE was not yet available. The present paper contributes to the analysis of the problem posed in [6, 7] as follows. By Theorem 4.1, the process $L_t := \exp(\alpha Y_t^B) > 0$ describes the maximal expected utility at time $t$ up to a deterministic factor, and our BSDE results ensure existence of a unique solution $Y^B$ to the BSDE (4.14). By Itô's formula, the BSDE for $Y^B$ [or (4.21)] implies that $L$ satisfies the BSDE

$$(4.29) \quad \begin{aligned} dL_t &= \frac{1}{2}\left(L_{t-}|\phi_t|^2 + \frac{|\widehat{\ell}_t|^2}{L_{t-}}\right) dt + \widehat{\ell} \, d\widehat{W} + \int_E \widetilde{\ell}_t \widetilde{\mu}(de, dt) \\ &= \frac{1}{2L_{t-}}|L_{t-}\phi_t + \widehat{\ell}_t|^2 + \widehat{\ell} \, dW + \int_E \widetilde{\ell}_t \widetilde{\mu}(de, dt), \end{aligned}$$

with $L_T = \exp(\alpha B)$, where $\widehat{\ell} := \alpha L_- Z^B$ and $\widetilde{\ell} := L_-(\exp(\alpha U^B(e)) - 1)$. Equation (4.29) corresponds to the key BSDEs (131) and (132) in [7], and our BSDE for $Y$ corresponds, up to multiplication by a constant, to their BSDE (136). The form of the BSDE (4.29) could also be motivated by the dual problem (4.20), using the fact that the density process of any $Q \in \mathbb{P}_f$ is a stochastic exponential driven by $W$ and $\widetilde{\mu}$ by (2.5) and following [6].

As already mentioned, the article [21] in some sense takes an opposite route by proving in a quite general model that existence and uniqueness for the particular BSDE corresponding to the exponential utility problem can be *derived from* the general duality results on the existence and structure of the solution to the exponential optimization problem in [8, 16].

4.3. *Dynamic utility indifference valuation and hedging.* In the same way as [21], we define the utility indifference value $\pi_t := \pi_t(B; \alpha)$ process for the claim $B$ under risk aversion $\alpha$ at any time $t \in [0, T]$ as the implicit solution of the equation

$$\text{(4.30)} \qquad V_t^{0,\alpha}(x) = V_t^{B,\alpha}(x + \pi_t), \qquad x \in \mathbb{R},$$

which relates the maximal expected utility functions of the optimization problems with and without ($B = 0$) a terminal liability. For exponential utility, it is clear from (4.12) that the solution $\pi$ to (4.30) does not depend on $x$. One should note that the notion "indifference value" is not uniform throughout the literature. In more classical terms, $\pi_t$ can be described as the offsetting variation of current wealth at time $t$ that compensates the investor for taking on the future liability $B$; see [10] for an exposition and references to the economic theory of value.

The utility indifference hedging strategy $\psi(B; \alpha)$ is defined as the difference of the respective optimal investment strategies,

$$\text{(4.31)} \qquad \psi := \psi(B; \alpha) := \theta^{B,\alpha} - \theta^{0,\alpha}.$$

Let $\widetilde{\mu}^E := \mu - \nu^E$ where

$$\text{(4.32)} \quad \nu^E(dt, de) := \exp(\alpha U^{0,\alpha}(e))\zeta(t,e)\lambda(de)\,dt =: \zeta^E(t,e)\lambda(de)\,dt$$

denotes the compensator of $\mu$ under the minimal entropy martingale measure $Q^E \equiv Q^{E,0}$, that is, $Q^{E,B}$ for $B = 0$. It it clear from (4.20) that $Q^E$ does not depend on $\alpha$. To see how this can be reconciled with (4.21), multiply by $\alpha$ the BSDE (4.14) from which $U^{0,\alpha}$ comes and use the uniqueness of the solution to conclude that $\alpha Y^{0,\alpha} = Y^{0,1}$, $\alpha Z^{0,\alpha} = Z^{0,1}$ and $\alpha U^{0,\alpha} = U^{0,1}$, then note that the $Q^{E,0}$-density (4.21) depends on $\alpha U^{0,\alpha}$.

For the remainder of this section, let us fix the claim $B \in L^\infty$ so we can ease the notation by omitting references to $B$ in some indices. The next theorem shows that the solution to the utility indifference pricing and hedging problem is characterized by the following single BSDE under the entropy minimal martingale measure $Q^E$:

$$\text{(4.33)} \quad \begin{aligned} Y_t &= B + \int_t^T \int_E \left( \frac{\exp(\alpha U_s(e))}{\alpha} - \frac{1}{\alpha} - U_s(e) \right) \zeta^E(t,e)\lambda(de)\,dt \\ &\quad - \int_t^T Z_s\,d\widehat{W}_s - \int_t^T \int_E U_s(e)\widetilde{\mu}^E(dt, de), \qquad t \in [0, T]. \end{aligned}$$

THEOREM 4.4. *The solution $\pi$ and $\psi$ to the dynamic utility indifference valuation and hedging problem for a claim $B \in L^\infty$ under risk aversion $\alpha > 0$ is described by the unique solution*

$$(Y^E, Z^E, U^E) := (Y^{E,\alpha}, Z^{E,\alpha}, U^{E,\alpha}) \in \mathcal{S}^\infty(Q^E) \times \mathcal{L}_T^2(Q^E) \times \mathcal{L}_\nu^2(Q^E)$$

to the BSDE (4.33) under $Q^E$. The utility indifference value process is

(4.34) $$\pi^\alpha = Y^{B,\alpha} - Y^{0,\alpha} = Y^{E,\alpha}$$

and the indifference hedging strategy is

(4.35) $$\psi^\alpha = Z^{B,\alpha} - Z^{0,\alpha} = Z^{E,\alpha},$$

where $Y^0, Y^B, Z^0$ and $Z^B$ are given by the BSDE solutions from Theorem 4.1 with terminal data $B$ and $0$, respectively.

At this point, a comment on a related result in [21] is instructive. Comparing the above BSDE to the one in equation (4.9) of [21] for a general filtration, one sees that the $\lambda$-integral part of the generator in our BSDE (4.33) corresponds to the compensator for the sum of jumps that appears in the BSDE in [21]. This sum is expressed directly in terms of jumps of $Y$ and therefore its compensator "makes it very hard to derive any properties." In our setting, on the other hand, the jump-related part of the generator can be explicitly expressed in terms of the integrand process $U$ and the $Q^E$-compensator $\nu^E$ of $\mu$; see (4.32). This form is highly amenable to further analysis as the subsequent results demonstrate.

PROOF OF THEOREM 4.4.　Let us first note that $\mathcal{L}^2_\nu(Q^E)$ equals $\mathcal{L}^2_{\nu^E}(Q^E)$, since the density $\exp(\alpha U^{0,\alpha})$ of $\nu^E$ with respect to $\nu$ is bounded from above and away from zero, and that $(\widehat{W}, \nu^E)$ with $Q^E$ fits into the setting of Section 3.

It follows directly from Theorem 4.1 and equations (4.30) and (4.31) that $\pi^\alpha = Y^B - Y^0$ and $\psi^\alpha = Z^B - Z^0$. It further follows from Theorem 4.1 that

$$(\delta Y, \delta Z, \delta U) := (Y^B - Y^0, Z^B - Z^0, U^B - U^0)$$

is the unique solution to the following BSDE (under $\widehat{P}$):

(4.36) $$\delta Y_t = B - \int_t^T \delta Z_s \, d\widehat{W}_s - \int_t^T \int_E \delta U_s(e) \widetilde{\mu}(ds, de) \\ + \int_t^T \int_E \bigg( \exp(\alpha U^0_s(e)) \frac{\exp(\alpha \delta U_s(e)) - 1}{\alpha} \\ - \delta U_s(e) \bigg) \zeta(s,e) \lambda(de) \, ds,$$

with $t \in [0, T]$. Recalling relation (4.32) and Theorem 12.28 from [13], a argument similar to that used for (4.25) yields that under a change of measure,

(4.37) $$\delta U \stackrel{\widehat{P}}{*} \widetilde{\mu} = \delta U \stackrel{Q^E}{*} \widetilde{\mu}^E + (\delta U(\exp(\alpha U^0) - 1)) * \nu.$$

Thereby, one can rewrite (4.36) to obtain a BSDE under $Q^E$ with a stochastic integral $\delta U * \widetilde{\mu}^E$. The remaining terms give rise to a different generator such that $(\delta Y, \delta Z, \delta U)$ is a solution to the BSDE (4.33) under $Q^E$, and by our BSDE results, such a solution is unique. $\square$

REMARK 4.5. As consequences of Theorem 4.4, we can observe several important and interesting properties of $\pi$ and $\psi$. For instance:

(1) The utility indifference value process $\pi$ is a $Q^E$-supermartingale since the generator in the BSDE (4.33) is nonnegative.

(2) The BSDE characterization of the solution implies that if $(\pi_t)$ and $(\psi_t)$ denote the solution from Theorem 4.4 with respect to claim $B$ at maturity $T$, then $(\pi_{t\wedge\tau})_{t\in[0,T]}$ and $(\psi_t \mathbb{1}_{[\![0,\tau]\!]}(t))_{t\in[0,T]}$ provide the indifference solution with respect to the claim $\pi_\tau$ at stopping time $\tau \leq T$. In this sense, the exponential utility indifference valuation and hedging approach is time consistent.

(3) Combining the results from Theorems 4.1 and 4.4, it is seen that $\pi$ and $\psi$ are the solution to a single optimization problem in the sense that

$$(4.38) \quad -\exp(\alpha \pi_t^\alpha) = \operatorname*{ess\,sup}_{\theta \in \Theta(Q^E, \alpha)} E_t^{Q^E}\left[ -\exp\left(-\alpha\left(\int_t^T \theta\, d\widehat{W} - B\right)\right)\right]$$

and the optimal strategy is attained by $\psi^\alpha \in \Theta(Q^E, \alpha)$, with $\Theta(Q^E, \alpha)$ being defined like $\Theta = \Theta(P, \alpha)$ from (4.8), but with $Q^E$ taking the role of $P$. This shows that the chosen definition of $\Theta$ transforms in a "good" way, and characterizes the solution $(\pi^\alpha, \psi^\alpha)$ to the utility indifference problem as the optimal solution to a single (primal) exponential utility optimization problem posed with respect to the minimal entropy measure $Q^E$ and over the set $\Theta(Q^E, \alpha)$ that is defined, consistently, with respect to $Q^E$. The latter, rather subtle, aspect of this statement appears to be new in the literature, while the general message from (4.38) corresponds to Proposition 3 in [21].

(4) The result of Theorem 4.4 implies, in combination with arguments from the proof of Theorem 4.2, that the utility indifference value process $\pi^\alpha$ is a (bounded) martingale under a suitable equivalent martingale measure $\widehat{Q}^B$ in $\mathbb{P}_e$ that depends on $B$ and $\alpha$. This confirms in our model framework an interesting observation made in [6, 7]; see the remark following their Proposition 27 in [7]. Thus far, it appears to be an open question to what extend such a property holds in general.

To prove it in our model, let $h(u) := \sum_{k=2}^\infty (\alpha u)^{k-1}/k!$. Then $h$ is a continuous function $\mathbb{R} \to \mathbb{R}$ with $h(u) > -1$ on $\mathbb{R}$ and $h(0) = 0$ such that $uh(u) = (\exp(\alpha u) - 1)/\alpha - u$, and the stochastic exponential $\mathcal{E}(h(U^E) * \widetilde{\mu}^E)$ is a martingale. The latter follows by Theorem 3.1 of [20] (or by Remark 3.1 in [17]) whose integrability condition is met because the compensator of

$((1+h(U^E))\log(1+h(U^E)) - h(U^E)) * \mu$ is bounded by the boundedness of $U$ and $\lambda(E)$. Hence,

$$d\widehat{Q}^B := \mathcal{E}(h(U^E) * \widetilde{\mu}^E)_T \, dQ^E$$

defines a probability measure. By the same Girsanov-type arguments as were used in the proof of Theorem 4.2, it follows that $\widehat{W}$ remains a Brownian motion under $\widehat{Q}^B$, while the compensator of $\mu$ under $\widehat{Q}^B$ becomes $\nu^{\widehat{Q}^B}(dt, de) = (h(U^E(e)) + 1)\nu^E(dt, de)$. By the same change-of-measure argument as used for (4.25) and (4.26), it follows from the BSDE (4.33) that $(Y^E, Z^E, U^E)$ is also the bounded solution to the following BSDE under $\widehat{Q}^B$:

$$(4.39) \qquad Y_t^E = B - \int_t^T Z_s^E \, d\widehat{W}_s + \int_t^T \int_E U_s^E(e) \widetilde{\mu}^B(ds, de),$$

where $\widetilde{\mu}^B := \mu - \nu^{\widehat{Q}^B}$ denotes the compensate measure under $\widehat{Q}^B$. Since $Y^E$ is bounded, it is clearly a martingale, and both stochastic integrals in (4.39) are $BMO(\widehat{Q}^B)$ martingales by Lemma 3.4.

Alternative proofs of some of the properties above, and others, in different or more general models, can be found in the literature; see [21] and the references therein.

4.4. *Asymptotics for vanishing risk aversion.* Finally, we prove that the utility indifference price and the indifference hedging strategy converge for vanishing risk aversion, in a suitable sense, to the conditional expectation process of $B$ and to the risk-minimizing strategy for $B$ under the minimal entropy martingale measure $Q^E$. This shows that such a convergence of the strategy, which has, to our best knowledge, thus far only been shown in [21] under the assumption of a continuous (Brownian) filtration, also holds in a setting like ours where the filtration is noncontinuous in that it allows for noncontinuous martingales. On this basis, one could conjecture that an asymptotic relation of this type should hold in general. We note that the convergence results of Section 4.4, that have been stated for vanishing risk aversion, could alternatively also be formulated for vanishing claim volume, by using a volume-scaling property for exponential utility (cf. Section 3 in [2]). This relates this article to a very interesting recent work of Kramkov and Sirbu [19] on the asymptotic of utility-based hedging strategies for small claim volumes. Unlike the present article, [19] does not investigate exponential utility but considers utility functions whose domain is the positive real half line.

Using suggestive notation anticipating Theorem 4.6, we denote by $(Y^{E,0}, Z^{E,0}, U^{E,0})$ from $\mathcal{S}^\infty(Q^E) \times \mathcal{L}_T^2(Q^E) \times \mathcal{L}_\nu^2(Q^E)$ the solution to the

BSDE (4.40) with zero generator under $Q^E$, with $t \in [0,T]$,

$$(4.40) \quad Y_t^{E,0} = B - \int_t^T Z_s^{E,0}\, d\widehat{W}_s - \int_t^T \int_E U_s^{E,0}(e)\widetilde{\mu}^E(ds, de).$$

It follows from this decomposition that $Z^{E,0}$ corresponds to the globally risk-minimizing strategy under the measure $Q^E$ (more precisely, to its risky asset's part) and that $Y^{E,0}$ is the associated valuation process (cf. [4]). For details on risk-minimization and relations to mean-variance hedging, see [23].

THEOREM 4.6. *Let* $(\pi^\alpha, \psi^\alpha, U^{E,\alpha}) = (Y^{E,\alpha}, Z^{E,\alpha}, U^{E,\alpha}) \in \mathcal{S}^\infty \times \mathcal{L}_T^2 \times \mathcal{L}_\nu^2$ *denote the BSDE solution (under $Q^E$) from Theorem 4.4 to the indifference pricing and hedging problem for risk aversion $\alpha \in (0, \infty)$. Then there is a constant $C = C(B) < \infty$ such that*

$$(4.41) \quad E_\tau^{Q^E}\Bigg[\sup_{u \in [\![\tau,T]\!]} |\pi_u^\alpha - Y_u^{E,0}|^2 + \int_\tau^T |\psi^\alpha - Z^{E,0}|^2\, ds \\ + \int_{]\!]\tau,T]\!] \times E} |U^{E,\alpha} - U^{E,0}|^2 \nu^E(ds, de)\Bigg] \leq \alpha^2 C,$$

*for all $\alpha \in (0,1]$ and all stopping times $\tau \leq T$. Hence,*

$$(4.42) \quad \begin{aligned} \sup_{t \in [0,T]} |\pi_t^\alpha - Y_t^{E,0}|^2 &\leq \alpha^2 C \qquad \text{for } \alpha \in (0,1] \\ \lim_{\alpha \downarrow 0} \int \psi^\alpha\, d\widehat{W} &= \int Z^{E,0}\, d\widehat{W} \qquad \text{in } BMO(Q^E), \\ \lim_{\alpha \downarrow 0} U^{E,\alpha} * \widetilde{\mu}^E &= U^{E,0} * \widetilde{\mu}^E \qquad \text{in } BMO(Q^E). \end{aligned}$$

*In particular,* $\lim_{\alpha \downarrow 0} \sup_{t \in [0,T]} |\pi_t^\alpha - Y_t^{E,0}| = 0$ *in* $L^\infty$, *and*

$$(4.43) \quad \begin{aligned} \lim_{\alpha \downarrow 0} \int \psi^\alpha\, d\widehat{W} &= \int Z^{E,0}\, d\widehat{W} \qquad \text{in } \mathcal{H}^2(Q^E), \\ \lim_{\alpha \downarrow 0} U^{E,\alpha} * \widetilde{\mu}^E &= U^{E,0} * \widetilde{\mu}^E \qquad \text{in } \mathcal{H}^2(Q^E). \end{aligned}$$

We note that the upper bound for $\alpha$ in $(0,1]$ is arbitrary in (4.41); the same result would hold for any finite bound other than 1, possibly with a different constant $C < \infty$.

PROOF OF THEOREM 4.6. Let $(Y^{E,\alpha}, Z^{E,\alpha}, U^{E,\alpha}) \in \mathcal{S}^\infty \times \mathcal{L}_T^2 \times \mathcal{L}_\nu^2$, for each $\alpha$, denote the BSDE solution from Theorem 4.4 which describes the solution to the indifference pricing and hedging problem. It follows, by applying (under $Q^E$) the estimate from Theorem 3.5 to the BSDE (4.33) and using Jensen's inequality, that there is a constant $c$ not depending on $\alpha$ such that

$$
\begin{aligned}
E_\tau^{Q^E} &\bigg[ \sup_{u \in [\![\tau,T]\!]} |Y_u^{E,\alpha} - Y_u^{E,0}|^2 \\
&+ \int_\tau^T |Z^{E,\alpha} - Z^{E,0}|^2\, ds + \int_{[\![\tau,T]\!] \times E} |U^{E,\alpha} - U^{E,0}|^2 \nu^E(ds, de) \bigg] \\
&\leq c E_\tau^{Q^E} \bigg[ \int_\tau^T \int_E \bigg| \frac{1}{\alpha} \exp(\alpha U_s^{E,\alpha}(e)) - \frac{1}{\alpha} - U_s^{E,\alpha}(e) \bigg|^2 \zeta^E(t,e) \lambda(de)\, dt \bigg]
\end{aligned}
$$

for all $\alpha \in (0,1]$. Since $|U^{E,\alpha}|$ is bounded uniformly in $\alpha$, for instance by $2b(0)$ from (3.12), the integrand $|(\exp(\alpha U_s^{E,\alpha}(e)) - 1)/\alpha - U_s^{E,\alpha}(e)|^2$ is bounded $P \otimes \nu$-a.e. by $\text{const} \cdot \alpha^2 < \infty$ for $\alpha \in (0,1]$, and as $\zeta^E$ is bounded and $\lambda(E)$ is finite, there is a constant $C < \infty$ such that the right-hand side of the above inequality can be bounded by $C\alpha^2$ for all $\alpha \in (0,1]$ and $\tau$. This yields (4.41)–(4.43), by letting $\tau = 0$ to obtain (4.43) and by using the characterization of BMO-martingales for (4.42); see Chapter 10 in [13]. □

4.5. *Examples.* This section outlines some areas where the general indifference results from Section 4 can be applied, and points out further connections to some closely related contributions in the literature.

One area of application involves incomplete stochastic volatility models. Consider, for instance, an increasing pure jump Lévy process $L$ without drift (i.e., a pure jump subordinator) and an independent Brownian motion $W$. Let $\mathbb{F} := \mathbb{F}^{(W,L)} = \mathbb{F}^{W+L}$ and let the price $S$ of the single risky asset evolve as

$$
\begin{aligned}
dS_t &= S_t \gamma_t(Y_{t-})\, dt + S_t \sigma_t(Y_{t-})\, dW_t, \qquad S_0 = s \in (0,\infty) \\
dY_t &= -KY_t + dL_t, \qquad Y_0 = y \in (0,\infty) \quad \text{and} \quad K \in (0,\infty),
\end{aligned}
\tag{4.44}
$$

for suitable $\mathcal{P} \otimes \mathcal{B}((0,\infty))$-measurable functions $\gamma, \sigma : \Omega \times [0,T] \times (0,\infty) \to \mathbb{R}$, with $\sigma > 0$. If $L$ is of finite activity and $|\frac{\gamma_t}{\sigma_t}(Y_{t-})|$ is bounded, this fits comfortably into the framework of Section 4.1 with $\mu := \mu^L$ and $\nu(dt, de) = \lambda(de)\, dt$ where $\lambda$ denotes the Lévy measure of $L$. Moreover, we can dispense with the independence assumption, for example, by a change to a new measure with $dP^{\text{new}}/dP^{\text{old}} = \mathcal{E}((\zeta-1) * (\mu - \nu))_T$ for a predictable bounded function $\zeta \geq 0$ such that the compensator of $\mu$ becomes $d\nu^{\text{new}} = \zeta\, d\nu^{\text{old}}$ and hence can depend on the history of $(W, L)$; see Example 2.1. In [5], the density process of the minimal entropy martingale measure $Q^E$ has been derived for a model like (4.44) with a specific Markovian choice of $\gamma$ and $\sigma$. Their results are both more and less general than ours. While they, for example, impose only a certain integrability on the market price of risk and do not require that $L$ is of finite activity, our framework can deal with other (also non-Markovian) choices of $\gamma$ and $\sigma$ and does not require the independence of the driving processes $W$ and $L$ under the objective measure. Moreover, we also obtain the density process of the optimal measure $Q^{E,B}$ for the dual problem if there is an additional liability $B \neq 0$.

Another application may be a "regime-switching" model

$$dS/S = \gamma(\eta_{t-})\,dt + \sigma(\eta_{t-})\,dW$$

where the local drift and volatility of $S$ are modulated by the state of a continuous finite-state Markov chain $(\eta_t)_{t \in [0,T]}$, independent of $W$, and where the claim $B$ could depend on the joint evolution of $(S, \eta)$. Under suitable assumptions on $\gamma, \sigma$, this fits well into the framework of Section 4 with $\mu := \mu^\eta$. Depending on the interpretation of $\eta$, say as an economic regime or a credit state, one arrives at an incomplete stochastic volatility model or at a model with tradable market risk and additional nontradable default risk. Again, the independence assumption between $\eta$ and $W$ could be considerably relaxed and mutual dependencies between $\eta$ and $S$ can be introduced; see [3] for a Markovian model of such type with utility indifference solutions in PDE-form that compare nicely to our BSDEs, and see [4] for similar results on risk-minimization that correspond to the limiting case for vanishing risk aversion from Section 4.4.

Another interesting area of application involves valuation and hedging problems for a portfolio of insurance policies (or of defaultable securities). To this end, it is common to consider a cumulative loss process $L_t = \sum_i \xi_i \mathbb{1}_{[\![T_i, T]\!]}(t)$ with random times $0 < T_1 < T_2 < \cdots$ where $T_i$ models the time of the $i$th insurance claim (or default event) and the random variable $\xi_i > 0$ may represent the size of claim $i$ (or the loss fraction times the notional). Suppose that $\mu := \mu^L$ has a compensator $\nu$ of the form (2.2)–(2.3) with respect to its usual filtration $\mathbb{F}^L$, and that $W$ is a Brownian motion independent of $L$. By Example 2.1(2), $\widetilde{\mu} = \mu - \nu$ and $W$ then have the representation property (2.4) with respect to $\mathbb{F} := \mathbb{F}^{(W,L)}$. This fits into our framework, and Section 4 describes the utility indifference solution for claims $B \in L^\infty$ that can depend on the joint history of $L$ and $S$ from (4.1), the latter process being driven by $W$. In the context of default risk, it should furthermore be interesting to also consider situations where the timing $(T_i)$ of losses and possibly also the loss sizes $(\xi_i)$ are stochastically related to the evolution of some common factors which can be hedged partially (via $S$ or $\widehat{W}$, resp.). To this end, dependencies between $W$ and $L$ ($\mu$) can again be introduced by a suitable change of measure (see Example 2.1) such that the compensator $\nu$ of $\mu = \mu^L$ can depend by its density $\zeta$ on the past evolution of the loss process $L$ and of the financial market prices $S$; see (2.2). This permits, for example, self-exciting features of $L$ ($\mu$). Our results can describe the utility indifference solution in a model with both default timing risk (via $T_i$) and additional, possibly nonpredictable, recovery risk (via $\xi_i$). In the latter aspect, this goes beyond the results of [3], even for the Markovian case.

## APPENDIX: CONDITIONAL FUBINI

This appendix provides a specific version of a conditional Fubini theorem. It may be folklore, but we have not found a reference elsewhere proving the specific result that we need.

LEMMA A.1. *Suppose $(X_t)$ is a $\mathcal{F} \otimes \mathcal{B}([0,T])$-measurable process with right continuous paths on a stochastic basis $(\Omega, \mathcal{F}, (\mathcal{F}_t)_{t \in [0,T]}, P)$ with $T < \infty$ such that $(X_t | t \in [0,T])$ is uniformly integrable. Let $u \in [0,T]$ and let $\sigma, \tau$ be stopping times with $\sigma \leq \tau \leq T$. Then there exists a measurable function*

$$F_\tau : (\Omega \times [0,T], \mathcal{F}_\tau \otimes \mathcal{B}([0,T])) \to (\bar{\mathbb{R}}, \mathcal{B}(\bar{\mathbb{R}}))$$

*such that $F_\tau(\cdot, t)$ is a version of $\mathbb{1}_{[\![\sigma \vee u, T]\!]}(t) E_\tau[X_t]$ for each $t$, and*

$$E_\tau \left[ \int_{\sigma \vee u}^T X_t \, dt \right] = \int_{\sigma \vee u}^T F_\tau(\cdot, t) \, dt = \int_{\sigma \vee u}^T E_\tau[X_t] \, dt.$$

PROOF. By suitable scaling, we may assume $T = 1$. We first show the claim for $\sigma = 0$ and $u = 0$. For $n \in \mathbb{N}$, define $F_\tau^n : \Omega \times [0,T] \to \mathbb{R}$, by

$$F_\tau^n(\omega, t) := E_\tau[X_{(j+1)/2^n}](\omega) \qquad \text{for } t \in \left[ \frac{j}{2^n}, \frac{j+1}{2^n} \right),$$

using the same version of $E_\tau[X_{(j+1)/2^n}]$ for all $t \in [j/2^n, j+1/2^n)$, and let $F_\tau^n(\omega, T) := E_\tau[X_T](\omega)$ for the terminal time $T = 1$. Then $F_\tau^n$ is $\mathcal{F}_\tau \otimes \mathcal{B}([0,T])$-measurable. Using right continuity of paths and uniform integrability, the conditional dominated convergence theorem yields that $F_\tau^n(t) = E_\tau[X_{q^n(t)}] \to E_\tau[X_t]$ a.s., since $q^n(t) := (\lfloor t 2^n \rfloor + 1)/2^n \searrow t$ for each $t \in [0,T)$ when $n \to \infty$. Hence, $F_\tau(\omega, t) := \limsup_n F_\tau^n(\omega, t)$ provides the required measurable version.

By uniform integrability of $X$, $E[\int_0^T |F_\tau(\cdot, t)| \, dt] \leq E[\int_0^T |X_t| \, dt]$ is finite. By Fubini's theorem, $t \mapsto F_\tau(t, \omega)$ is for a.a. $\omega$ integrable with respect to the Lebesgue measure, and $\int_0^T F_\tau(t, \omega) \, dt$ is $\mathcal{F}_\tau$-measurable. Finally, again by Fubini's theorem, we obtain

$$(A.1) \quad E\left[ I_A \int_0^T X_t \, dt \right] = \int_0^T E[\mathbb{1}_A F_\tau(\cdot, t)] \, dt = E\left[ \mathbb{1}_A \int_0^T F_\tau(\cdot, t) \, dt \right]$$

for any $A \in \mathcal{F}_\tau$, which yields the claim for the case $\sigma, u = 0$. With $F_\tau$ from above, we construct the solution for $u \in [0,T]$ and $\sigma \leq \tau$ as follows. Let $\widetilde{X} := \mathbb{1}_{[\![\sigma \vee u, T]\!]} X$, and define $\widetilde{F}_\tau(\omega, t) := \mathbb{1}_{[\![\sigma \vee u, T]\!]}(t) F_\tau(\omega, t)$. Noting that $\mathbb{1}_{[\![\sigma \vee u, T]\!]}(t) = \mathbb{1}_{[\![\sigma,T]\!]}(t)\mathbb{1}_{[u,t]}(t) \in \mathcal{F}_\sigma \subset \mathcal{F}_\tau$, one sees that $\widetilde{F}$ is $\mathcal{F}_\tau \otimes \mathcal{B}([0,T])$-measurable, $\widetilde{F}_\tau(\cdot, t)$ is a version of $\mathbb{1}_{[\![\sigma \vee u, T]\!]}(t)E_\tau[X_t] = E_\tau[\widetilde{X}_t]$ for all $t$ and a Fubini argument like (A.1), but with $\widetilde{X}$ and $\widetilde{F}_\tau$, completes the proof. $\square$

DEPARTMENT OF MATHEMATICS
IMPERIAL COLLEGE
LONDON SW7 2AZ
UNITED KINGDOM
E-MAIL: dirk.becherer@ic.ac.uk